\documentclass[11pt]{article}
\usepackage[latin9]{inputenc}
\usepackage{geometry}
\usepackage{amsmath}
\usepackage{amssymb}
\usepackage{algorithm}
\usepackage{algorithmic}
\usepackage{float}
\usepackage{graphicx}
\usepackage{mathtools}
\usepackage{pgffor}
\usepackage[numbers, sort]{natbib}

\makeatletter

\usepackage{amsfonts}\usepackage{nopageno}
\newcommand{\diag}{{\rm diag}}
\reversemarginpar

\makeatother
\let\clearpage\relax

\begin{document}
\title{A Vectorial Approach to Unbalanced Optimal Mass Transport}
\author{Jiening Zhu, Rena Elkin, Jung Hun Oh, Joseph O. Deasy, Allen Tannenbaum
\thanks{J.\ Zhu is with the Department of Applied Mathematics \& Statistics, Stony Brook University, NY; email: jiening.zhu@stonybrook.edu}
\thanks{R. Elkin, J-H Oh, and J. Deasy are with the Department of Medical Physics, Memorial Sloan Kettering Cancer Center, NY; email: ElkinR@mskcc.org, OhJ@mskcc.org, DeasyJ@mskcc.org}
\thanks{A.\ Tannenbaum is with the Departments of Computer Science and Applied Mathematics \& Statistics, Stony Brook University, NY; email: allen.tannenbaum@stonybrook.edu}}
\date{\today}
\maketitle
\begin{abstract}
Unbalanced optimal mass transport (OMT) seeks to remove the conservation of mass constraint by adding a source term to the standard continuity equation in the Benamou-Brenier formulation of OMT. In this note, we show how the addition of the source fits into the vector-valued OMT framework.
\end{abstract}

\section{Introduction}
\par
Optimal mass transport (OMT) is a very important subject in mathematics, originating with the French civil engineer and mathematician Gaspard Monge in 1781 \cite{villani1, villani2}. Recently, the theory has undergone a massive growth, with numerous applications to various areas including signal processing, machine learning, computer vision, meteorology, statistical physics, quantum mechanics, and network theory \cite{Arjovsky2017,Haker2004,rr,inbook,Statement2020}. Several works deal with extensions of the theory to the  unbalanced and vector-valued cases; see \cite{Chizat, vectorvalued} and the references therein. In standard treatments, vector-valued and unbalanced extensions are treated separately. In the present work, we show that unbalanced OMT can fit into the model of vector OMT by taking a special set of weight parameters, which gives us fresh way to treat the problem.

In the $L^2$ setting, both extensions arise from the computational fluid dynamics (CFD) approach to OMT by Benamou and Brenier \cite{BB}, which was a major development in optimal mass transport theory. For ``square of distance'' based cost functionals, the Benamou and Brenier framework is equivalent to the original one while giving an explicit interpolation of two mass distributions regarded as a geodesic on the space of probability distributions \cite{Otto}. The CFD method views the optimal mass transport as the minimizer of a kinetic energy functional subject to a continuity constraint. By modifying the energy functional and the constraint, unbalanced OMT and vector OMT versions were developed \cite{chizat:hal-01271981,Chizat,vectorvalued}. More specifically, by adding a source term, unbalanced OMT can handle transport problems when the mass is not conserved. That is, mass can be created or destroyed during the process. By introducing the divergence of the flow between channels, vector OMT extends the density from scalar to vector-valued and even matrix-valued \cite{CheGeoTan16,cheGeoTan16b}.

In this note, we indicate how to reformulate unbalanced OMT into a vector-valued OMT framework, and thus may be implemented through any existing code for vector OMT. All that needs to be done is adding a new weight matrix without changing the main structure of the code.

In this note, we focus on the connection between these two extensions. We will first give the background on the two extensions of OMT. After that we will show how to reformulate unbalanced OMT to a vector OMT problem. We will also look into reformulating unbalanced vector OMT to a common vector OMT problem. We conclude with some illustrative numerical results.

\section{Background}
In this section, we briefly introduce the Benamou and Brenier approach of OMT \cite{BB} as well as the two extensions mentioned above.

\subsection{Benamou and Brenier}
The original formulation of OMT due to Monge \cite{villani1, villani2} may be expressed as follows:
\begin{equation}
\label{1}
  \inf_T\{\int_{E}c(x,T(x))\rho_0(x)dx\ |\ T_{\#}\rho_0=\rho_1\},
\end{equation}
where $c(x,y)$ is the cost of moving unit mass from $x$ to $y$, $T$ is the transport map, and $\rho_0,\rho_1$ are two probability distributions defined on $E$ a subdomain of $\mathbb{R}^n$. $T_\#$ denotes the push forward of $T$.

As pioneered by Leonid Kantorovich \cite{villani1, villani2}, the Monge formulation of OMT may be relaxed replacing transport maps $T$ by couplings $\pi$:
\begin{equation}
  \inf_{\pi\in\Pi(\rho_0,\rho_1)}\int_{E}c(x,y)\pi(dx,dy),
\end{equation}
where $\Pi(\rho_0,\rho_1)$ denotes the set of all the couplings between $\rho_0$ and $\rho_1$ (joint distributions whose two marginal distributions are $\rho_0$ and $\rho_1$).
One may show that for $c(x,y) =\|x-y\|^2$ (square of distance function), the Kantorovich and Monge formulations are equivalent; see \cite{villani1, villani2} and the references therein.

Moreover for $c(x,y)=||x-y||^2$, the specific infimum is called Wasserstein 2 distance ($\mathcal{W}_2$). Benamou and Brenier pointed out that this can be written in an equivalent computational fluid dynamic formulation:
\begin{subequations}
  \begin{align}
  \mathcal{W}_2(\rho_0,\rho_1)^2=&\inf_{\rho,v}\int_{0}^{1}\int_{E}\rho(t,x)||v(t,x)||^2dxdt\\
  &\frac{\partial \rho}{\partial t}+\nabla\cdot(\rho v)=0\\
  & \rho(0,\cdot)=\rho_0,\rho(1,\cdot)=\rho_1.
\end{align}
\end{subequations}
Thus, the latter reformulates OMT as one of minimizing a kinetic energy functional subject to a continuity constraint. The continuity constraint states that the change of density at each point over time is due to the flux of its neighborhood.

The Benamou-Brenier optimal solution may be related to the the original one of Monge. Indeed, if we have the optimal transport plan $T$ in (\ref{1}), the interpolation can be expressed as $\rho(t,\cdot)=((t\cdot T+(1-t)\cdot id)_{\#}\rho_0$.
For computational purposes, the Benamou-Brenier model is usually written as a convex optimization problem in momentum form ($p=\rho v$):
\begin{subequations}
  \begin{align}
  \mathcal{W}_2(\rho_0,\rho_1)^2=&\inf_{\rho,p}\int_{0}^{1}\int_{E}\frac{p(t,x)^2}{\rho(t,x)}dxdt\\
  &\frac{\partial \rho}{\partial t}+\nabla\cdot p=0\\
  & \rho(0,\cdot)=\rho_0,\rho(1,\cdot)=\rho_1.
\end{align}
\end{subequations}

\subsection{Unbalanced OMT}

In the standard formulation, $\rho_0,\rho_1$ need to be balanced, i.e., have the same total mass. We enforced this by simply taking them to be probability distributions ($\int_{E}\rho_0(x)dx=\int_{E}\rho_1(x)dx=1$). In the unbalanced case, total mass is not required to be preserved. This is very important for a number of applications in which mass may be created or destroyed.

There are many ways to extend the original setting to the unbalanced case \cite{Gangbo2019}. In the present work, we will just treat the case in which a source term is
added under the Fisher-Rao smooth formulation \cite{Chizat}. Namely,
\begin{subequations}\label{unbalanced}
  \begin{align}
  \mathcal{W}_{FR}(\rho_0,\rho_1)^2=&\inf_{\rho,p}\int_{0}^{1}\int_{E}\frac{p(t,x)^2}{\rho(t,x)}+\gamma \frac{s(t,x)^2}{\rho(t,x)} dxdt\\
  &\frac{\partial \rho}{\partial t}+\nabla\cdot p=s\label{unbalanced_con}\\
  & \rho(0,\cdot)=\rho_0,\rho(1,\cdot)=\rho_1.
\end{align}
\end{subequations}

There is an extra source term $s$ in the continuity equation. The change of density at each point over time now is not only due to the flux of its neighborhood, but also a source. Mass can be created or destroyed at any point and time. The parameter $\gamma$ is a weight to control how much source one wants to use.

\subsection{Vector-valued OMT}

Vector-valued OMT considers vector-valued distributions instead of only the classical scalar ones \cite{vectorvalued}. With more than one channel, the change of mass is not only due to the flux of its neighborhood but also the redistribution within the given channels. Namely, we have
\begin{subequations}\label{vector1}
  \begin{align}
  \mathcal{W}_{V}(\rho_0,\rho_1)^2=& \inf_{\rho,p,u}\int_{0}^{1}\int_{E} p^{T}\diag(\rho)^{-1}p+\gamma u^T[\diag(\mathbb{F}^T_2\rho)^{-1}+\diag(\mathbb{F}^T_1\rho)^{-1}]u\ dx dt \\
  & \frac{\partial \rho}{\partial t}+\nabla_x\cdot p-\nabla_\mathcal{F}^*u=0\label{vec_con}\\
  &\rho(0,\cdot)=\rho_0,\quad \rho(1,\cdot)=\rho_1.
\end{align}
\end{subequations}

As mentioned above, there are two kinds of flows for the vector-valued density $\rho$. The first is the spatial flux $p$ and its corresponding divergence is denoted as $\nabla_x\cdot$, and the second is the flux $u$ between channels, whose  discrete divergence is denoted by $\nabla_\mathcal{F}^*$. $p$ and $u$ are both vectors. We note that $\mathbb{F}_1$ is the all-1 part of the incident matrix $\mathbb{F}$ (sources of all the edges) and $\mathbb{F}_2=\mathbb{F}_1-\mathbb{F}$ (sinks of all the edges). As $u$ is defined on every edge and the density is only defined on nodes, $\diag(\mathbb{F}^T_2\rho)^{-1}+\diag(\mathbb{F}^T_1\rho)^{-1}$ assigns a density to each edge using the density of two end points so that we can compute the kinetic energy for the flows on edges.

Vector-valued OMT can be considered as a general OMT problem on $E \times F$, where $E$ is the space on which the vector-valued density is defined and $F$ is a connected graph. With this understanding, we can rewrite the energy functional in the following equivalent form:
\begin{equation}\label{vector2}
    \mathcal{W}_{V}(\rho_0,\rho_1)^2=\inf_{\rho,p,u}\int_{0}^{1}\int_{E}\sum_{c\in Nodes(F)}\frac{p(t,x,c)^2}{\rho(t,x,c)}dx dt+\gamma\int_{0}^{1}\int_{E}\sum_{e\in Edges(F)}\frac{u(t,x,e)^2}{\tilde{\rho}(t,x,e)}dx dt,
\end{equation}
where $\tilde{\rho}$ is density of edge.

\subsection{Unbalanced vector-valued OMT}

In this section, we write down the unbalanced vector-valued version of OMT. Namely,
\begin{subequations}
  \begin{align}
  \mathcal{W}_{VS}(\rho_0,\rho_1)^2=& \inf_{\rho,p,u}\int_{0}^{1}\int_{E} p^{T}\diag(\rho)^{-1}p+\gamma u^T[\diag(\mathbb{F}^T_2\rho)^{-1}+\diag(\mathbb{F}^T_1\rho)^{-1}]u+\eta s^{T}\diag(\rho)^{-1}s \ dx dt\label{VS} \\
  & \frac{\partial \rho}{\partial t}+\nabla_x\cdot p-\nabla_\mathcal{F}^*u=s\\
  &\rho(0,\cdot)=\rho_0,\quad \rho(1,\cdot)=\rho_1.
\end{align}
\end{subequations}

Note that the source term is also a vector. So we need to write the energy functional for this vector term. Regarding the continuity equation, there are three possibilities for the change of mass over time: the flux of its neighborhood, flows between channels and source.

\section{Reformulation unbalanced scalar OMT to vector-valued OMT}

When rewriting vector OMT in (\ref{vector1}), we should note the similarity with unbalanced OMT. The flow between channels may be utilized as the source. The source term describes the creation or vanishing of mass. This can be modeled via an extra layer: the mass created originates from that layer and the vanishing mass just goes to that layer.

\subsection{Source layer: scalar case}

Here the input source and target are regarded as two mass distributions on the subdomain $E$. Total mass may or may not be preserved. We add an extra source layer which is parallel to original space. For each point in $E$, there is an edge connecting with source layer. As now we have an extra layer, we can put the difference of mass into the new layer so that the total mass of the 2-vector new structure is preserved. The difference of mass can just be distributed uniformly to source layer or put to some area of interest for further specific applications.

Now we view the flow between channels as source term. Consider the continuity equations of (\ref{unbalanced_con}) and (\ref{vec_con}):
\begin{align*}
  \frac{\partial \rho}{\partial t}+\nabla\cdot p=&s\\
  \frac{\partial \rho}{\partial t}+\nabla_x\cdot p=&\nabla_\mathcal{F}^*u.
\end{align*}
We can use $\nabla_\mathcal{F}^*u$ as the source term $s$. Because of the special graph structure (there is only one edge), we have:
\begin{equation}
  \nabla_\mathcal{F}^*u=u=s.
\end{equation}
Now consider the second term in (\ref{vector2}). Since there is only one edge, we can just omit the summation:
\begin{equation}
  \gamma\int_{0}^{1}\int_{E}\sum_{e\in Edges(F)}\frac{u(t,x,e)^2}{\tilde{\rho}(t,x,e)}dx dt=\gamma\int_{0}^{1}\int_{E}\frac{u(t,x)^2}{\tilde{\rho}(t,x)}dx dt.
\end{equation}
If we further take the density of each edge $\tilde{\rho}(t,x)=\rho(t,x)$, this term is exactly the same as the second term in the energy functional in unbalanced OMT setting.

Now we consider the first term of the energy functional:
\begin{equation}
  \int_{0}^{1}\int_{E}\sum_{c\in Nodes(F)}\frac{p(t,x,c)^2}{\rho(t,x,c)}dx dt=\int_{0}^{1}\int_{E}\frac{p(t,x,c_1)^2}{\rho(t,x,c_1)}+\frac{p(t,x,c_2)^2}{\rho(t,x,c_2)}dx dt.
\end{equation}
There are only two channels, i.e., $c_1$ denotes the original channel and $c_2$ denotes the source layer. Notice that the integral with respect to $c_1$ term is exactly the first term of (\ref{unbalanced}). We want to get rid of the $c_2$ integral. This can be done by introducing a small weight parameter for that layer.

\subsection{A naive generalization of vector-valued OMT: addition of  weight}

In the original setting, one treats the kinetic energy of different layers in an identical manner. But by simply introducing a weight, we can treat each layer differently:
\begin{equation}\label{vector3}
    \mathcal{W}_{V'}(\rho_0,\rho_1)^2=\inf_{\rho,p,u}\int_{0}^{1}\int_{E}\sum_{c\in Nodes(F)}w(c)\frac{p(t,x,c)^2}{\rho(t,x,c)}dx dt+\gamma\int_{0}^{1}\int_{E}\sum_{e\in Edges(F)}\frac{u(t,x,e)^2}{\tilde{\rho}(t,x,e)}dx dt,
\end{equation}
where $w(c)$ is the weighting parameters. We can choose different values for different channels. For another form of vector OMT:
\begin{equation}\label{vector4}
  \mathcal{W}_{V'}(\rho_0,\rho_1)^2=\inf_{\rho,p,u}\int_{0}^{1}\int_{E} p^{T}\diag(w)\diag(\rho)^{-1}p+\gamma u^T[\diag(\mathbb{F}^T_2\rho)^{-1}+\diag(\mathbb{F}^T_1\rho)^{-1}]u\ dx dt.
\end{equation}

Under this setting, if we take $w(c_1)=1$ and $w(c_2)$ very small, then the integral in the source layer is very small in that energy functional. Clearly, unbalanced OMT is almost equivalent to the extra source layer vector OMT with the latter set of weight parameters.

\subsection{Implementation: scalar case}

\begin{figure}[H]
  \centering
  \includegraphics[width=0.6\linewidth]{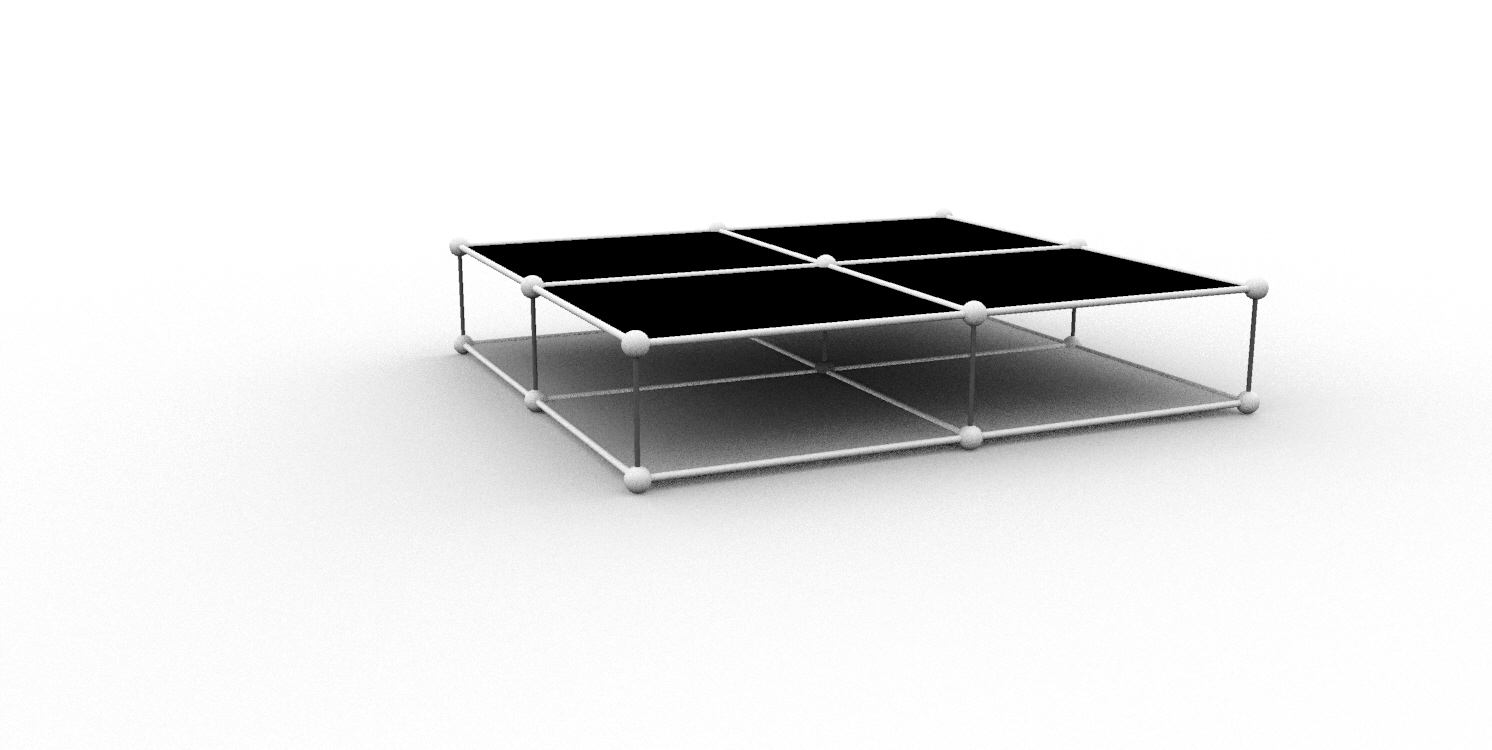}
  \caption{The bottom layer (gray) is the source layer, the weight for the flow within it is very small. The flow on the edge between two layers is our source.} \label{fig:source}
  \end{figure}

By introducing the source layer, there is an edge between two channels for which we use the flow on that edge as source,  and thus there is an extra integral on that new layer. We add very small weight parameters so that two energy functionals are almost identical. See Figure~\ref{fig:source}.

It is quite straightforward to implement unbalanced OMT from vector-valued OMT code. We need to only alter a few lines of code to make it work for unbalanced case:
\begin{enumerate}
  \item Initialization:
  From the input, first construct two extended structures for vector OMT. The first layer is the original input and the second layer contains the mass difference. If the starting density distribution has more total mass, then the mass difference is added to the second layer of target density distribution. If the starting density distribution has less total mass, then the mass difference is added to the second layer of starting density distribution itself.
  \item Set weighting parameters:
  Employing the code for the energy functional, just add weighting parameters to the corresponding layers.
\end{enumerate}

\section{Reformulation of unbalanced vector OMT to vector-valued OMT}
The reformulation is very similar to the scalar case. The same simple idea is to add a new source layer which connects to each of the original layers. The flows on those added edges are our sources.

\subsection{Source layer: vector-valued case}

In the vector-valued case, the input source and target are two $n$-vector valued distributions. Total mass may or may not be preserved. We add an extra source layer which is connected to each of the existing layers. As before, we can put the difference of mass into the new layer so that the total mass of the $(n+1)$-vector new structures is preserved.

We call the incident matrix for those newly added edges $\mathcal{G}$, which is a submatrix of the new large incident matrix $\tilde{\mathcal{F}}$. We just split $\tilde{\mathcal{F}}$ for illustration purpose. They can be easily pieced together as in original vector setting:
\begin{equation}
  \nabla_{\tilde{\mathcal{F}}}^*\left[\begin{array}{c}
                                      u \\
                                      v
                                    \end{array}\right]=\left[\nabla_\mathcal{F}^*\ \nabla_\mathcal{G}^*\right]\left[\begin{array}{c}
                                      u \\
                                      v
                                    \end{array}\right]=\nabla_\mathcal{F}^*u+\nabla_\mathcal{G}^*v,
\end{equation}
where $u$ is the flow within the existing edges and $v$ is the flow within the newly added edges.

Now we consider the continuity equations of unbalanced vector OMT and vector OMT with a new layer:
\begin{align*}
  \frac{\partial \rho}{\partial t}+\nabla\cdot p-\nabla_\mathcal{F}^*u=&s \\
  \frac{\partial \rho}{\partial t}+\nabla_x\cdot p-\nabla_\mathcal{F}^*u=&\nabla_\mathcal{G}^*v.
\end{align*}
We can use $\nabla_\mathcal{G}^*v$ as the source $s$. Now similar to the scalar case:
\begin{equation}\label{relation}
  \nabla_\mathcal{G}^*v=v=s.
\end{equation}
The energy functional now becomes:
\begin{align}
\mathcal{W}_{V''}(\rho_0,\rho_1)^2=\inf_{\rho,p,u,v}\int_{0}^{1}\int_{E}\sum_{c\in Nodes(\tilde{\mathcal{F}})}\frac{p(t,x,c)^2}{\rho(t,x,c)}dx dt+&\gamma\int_{0}^{1}\int_{E}\sum_{e\in Edges(\mathcal{F})}\frac{u(t,x,e)^2}{\tilde{\rho}(t,x,e)}dx dt\\
+&\eta\int_{0}^{1}\int_{E}\sum_{e\in Edges(\mathcal{G})}\frac{v(t,x,e)^2}{\tilde{\rho}(t,x,e)}dx dt.\notag
\end{align}
With the above relationship (\ref{relation}), it is easy to see that the last term exactly fits the original source energy term (\ref{VS}).

Again there is a new part of kinetic energy due to the newly added source layer. We can introduce weight parameters to make the extra term almost disappear. Same as adding weight technique for different layers, we can put $\gamma$ and $\eta$ together as a weight matrix for different edges so that the final form is more concise:
\begin{subequations}
  \begin{align}
\mathcal{W}_{V''}(\rho_0,\rho_1)^2=&\inf_{\rho,p,\tilde{u}}\int_{0}^{1}\int_{E} p^{T}\diag(w_1)\diag(\rho)^{-1}p+ \tilde{u}^T \diag(w_2)\diag(\tilde{\rho})^{-1}\tilde{u}\ dx dt\\
  & \frac{\partial \rho}{\partial t}+\nabla_x\cdot p-\nabla_{\tilde{\mathcal{F}}}^*\tilde{u}=0\\
  &\rho(0,\cdot)=\rho_0,\quad \rho(1,\cdot)=\rho_1,
\end{align}
\end{subequations}
where $\diag(w_1)$ denotes the weighting matrix for different layers and $\diag(w_2)$ denotes the weighting matrix for the edges we add (weight $=\eta$) and existing edges (weight$=\gamma$). $\tilde{u}=\left[\begin{array}{c}
                                      u \\
                                      v
                                    \end{array}\right]$ is all the flows that between channels within the new large graph $\tilde{\mathcal{F}}$ and $\tilde{\rho}$ is the density assigned to each edge.

\subsection{Implementation: vector-valued case}

\begin{figure}[H]
  \centering
  \includegraphics[width=0.6\linewidth]{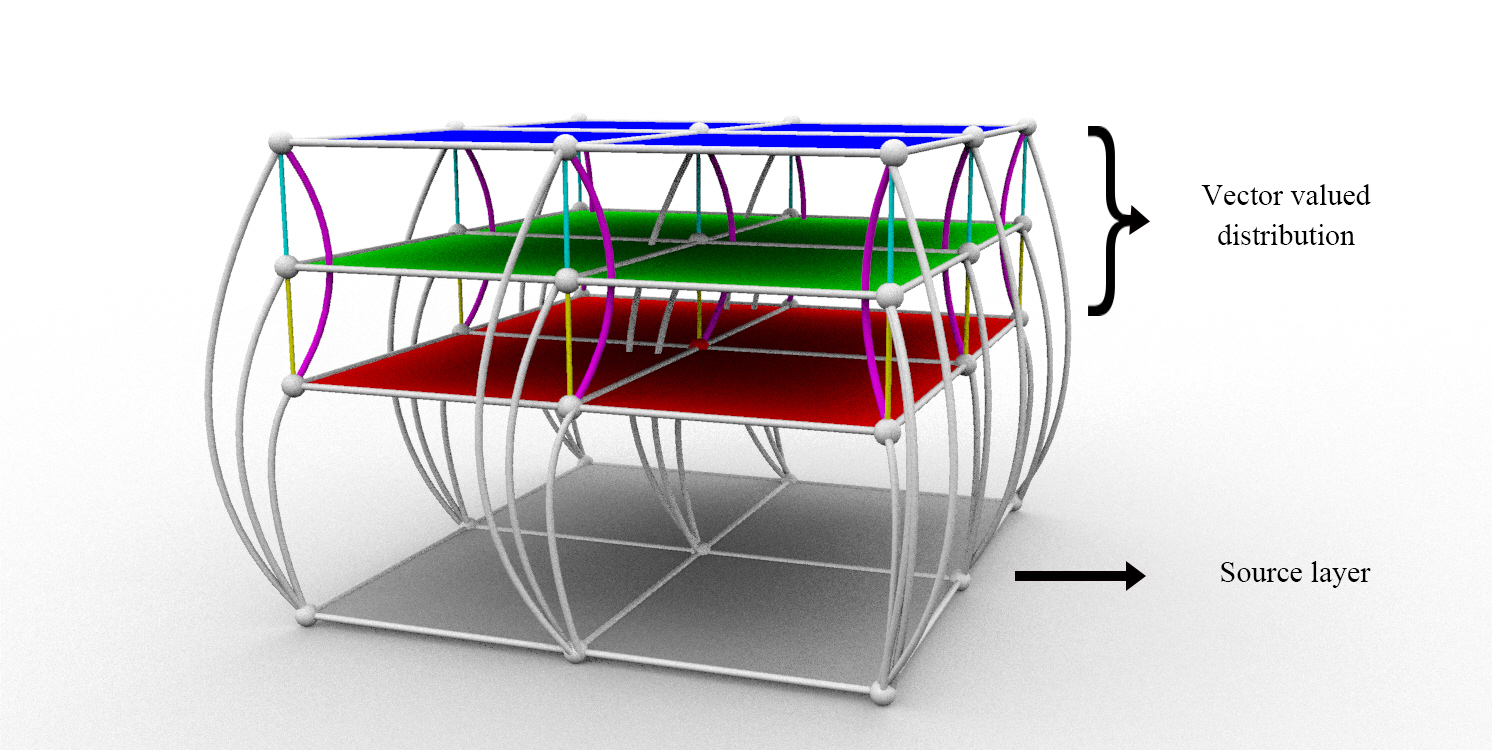}
  \caption{The bottom layer (gray) is the source layer, the weight for the flow within it is very small. The flow on the edge between the source layer and other existing layers are our sources.} \label{fig:figure2}
\end{figure}

By introducing a source layer, there is an edge between the source layer and each existing channel. We use the flows on those edges as sources. As before, we add very small weight parameters for the source layer so that two energy functionals are almost identical. See Figure~\ref{fig:figure2}.

As before, it is also very easy to implement unbalanced vector OMT from vector OMT code.
\begin{enumerate}
  \item Initialization:
  From the input, first construct the extended structures for vector OMT. Add an extra layer. Connect that new layer with each of other layers. Put the difference of total mass into that newly added layer.
  \item Set weight parameters:
  Employing the code for the energy functional, add weighting parameters to corresponding layer and corresponding edges (existing edges and newly added edges).
\end{enumerate}
With the simple change of the original vector-valued  OMT, we can use the same code for unbalanced vector OMT case without even touching the main structure of the original code.

\section{Numerical results}
We tested our new formulation on several images using the numerical algorithm from \cite{CheKaoHabGeoTan17}. Gray scale images are general mass distributions on a rectangular area while color images are vector valued distributions.

\subsection{Unbalanced OMT}
We tested an example of moving two Gaussian distributions. Though this example has preserved total mass, the optimal solution of OMT still uses the source term. We can look at the effect of tuning the weight parameter $\gamma$.
The source image and target images are:
\begin{figure}[H]
  \centering
  \includegraphics[width=0.2\linewidth]{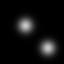}
  \includegraphics[width=0.2\linewidth]{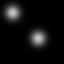}
  \caption{Source and target distributions}
\end{figure}
With large $\gamma$, it is expensive to use the source term, so the source layer is almost zero at all times.
\begin{figure}[H]
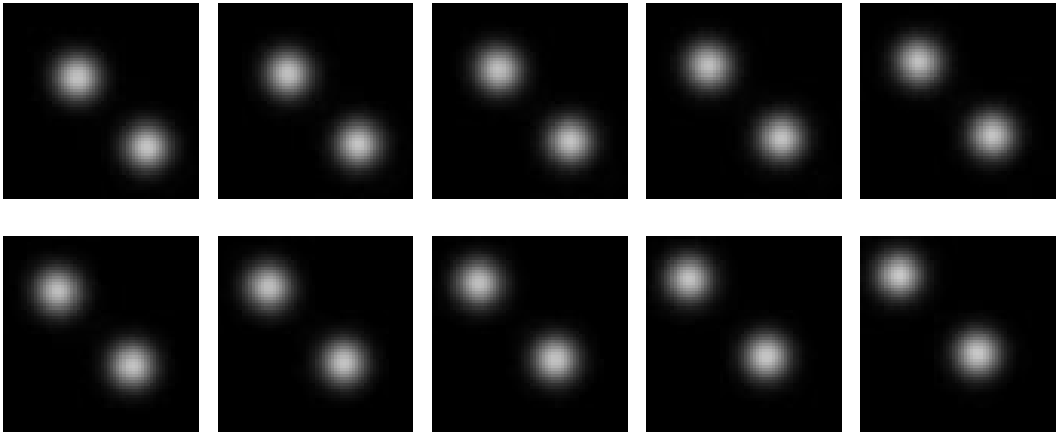

  \foreach \t in {1,2,3,4,5}{
  \includegraphics[width=0.17\linewidth]{figures//01_\t.jpg}
  }\\ \\
  \foreach \t in {6,7,8,9,10}{
  \includegraphics[width=0.17\linewidth]{figures//01_\t.jpg}
  }
  \caption{Density over time}
\end{figure}
\begin{figure}[H]
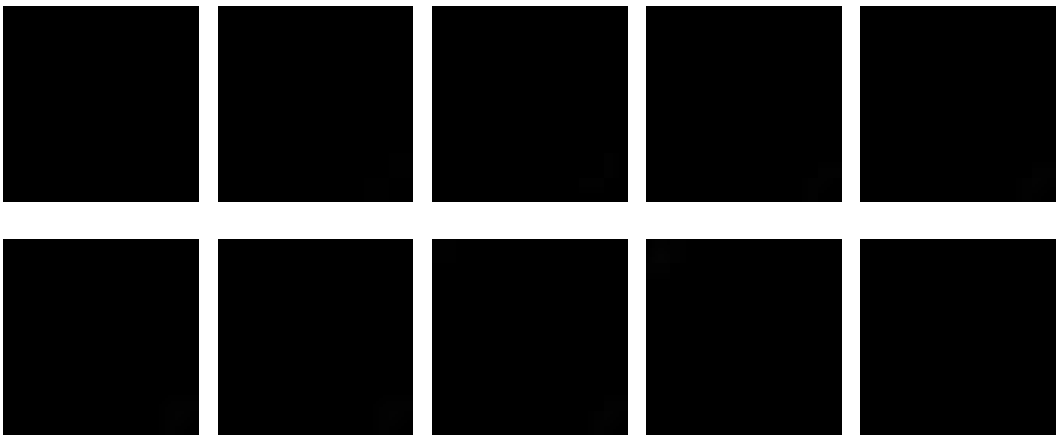

  \foreach \t in {1,2,3,4,5}{
  \includegraphics[width=0.17\linewidth]{figures//01_s_\t.jpg}
  }\\ \\
   \foreach \t in {6,7,8,9,10}{
  \includegraphics[width=0.17\linewidth]{figures//01_s_\t.jpg}
  }
  \caption{Source over time}
\end{figure}
With small $\gamma$, it is cheap to use the source term, thus much of the mass goes through the source layer.
\begin{figure}[H]
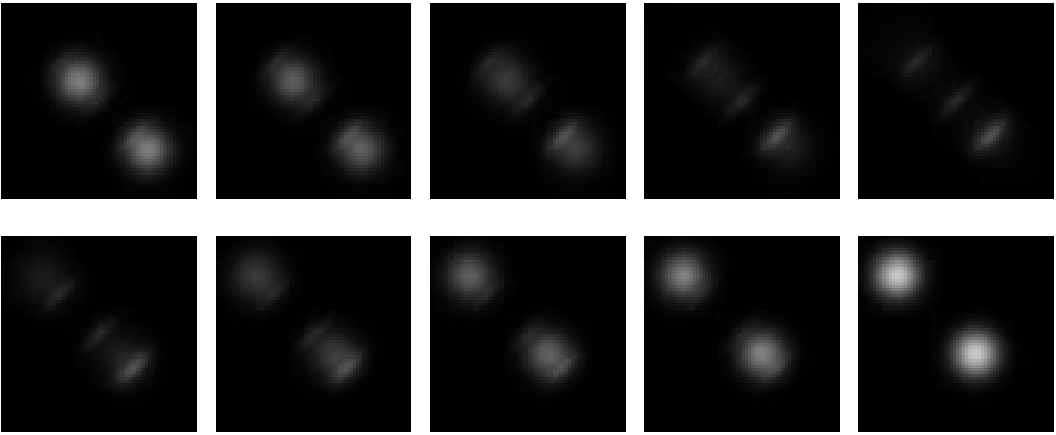

  \foreach \t in {1,2,3,4,5}{
  \includegraphics[width=0.17\linewidth]{figures//00001_\t.jpg}
  }\\ \\
  \foreach \t in {6,7,8,9,10}{
  \includegraphics[width=0.17\linewidth]{figures//00001_\t.jpg}
  }
  \caption{Density over time}
\end{figure}
\begin{figure}[H]
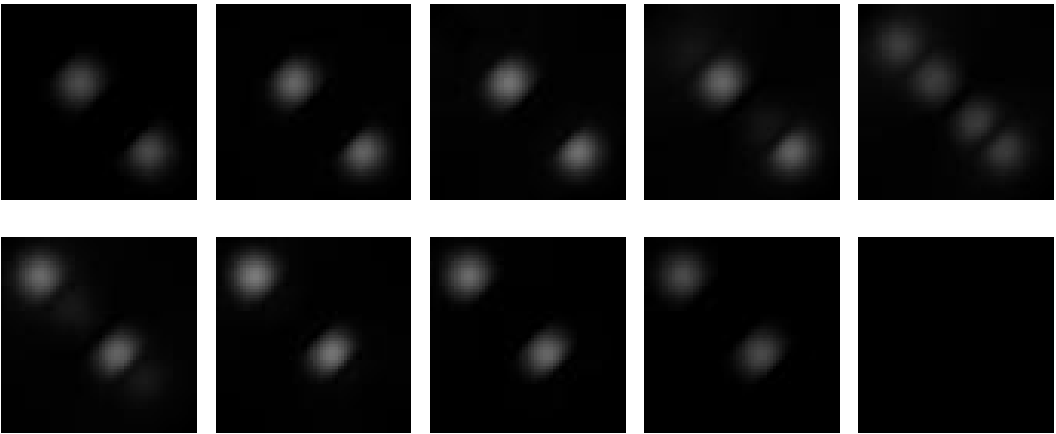

  \foreach \t in {1,2,3,4,5}{
  \includegraphics[width=0.17\linewidth]{figures//00001_s_\t.jpg}
  }\\ \\
  \foreach \t in {6,7,8,9,10}{
  \includegraphics[width=0.17\linewidth]{figures//00001_s_\t.jpg}
  }
  \caption{Source over time}
\end{figure}
We also tested two images with very large total mass difference:
\begin{figure}[H]
  \centering
  \includegraphics[width=0.2\linewidth]{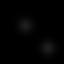}
  \includegraphics[width=0.2\linewidth]{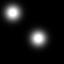}
  \caption{Source and target distributions}
\end{figure}
\begin{figure}[H]
  \foreach \t in {1,2,3,4,5}{
  \includegraphics[width=0.17\linewidth]{figures//l_\t.jpg}
  }\\ \\
  \foreach \t in {6,7,8,9,10}{
  \includegraphics[width=0.17\linewidth]{figures//l_\t.jpg}
  }
  \caption{Density over time}
\end{figure}
\begin{figure}[H]
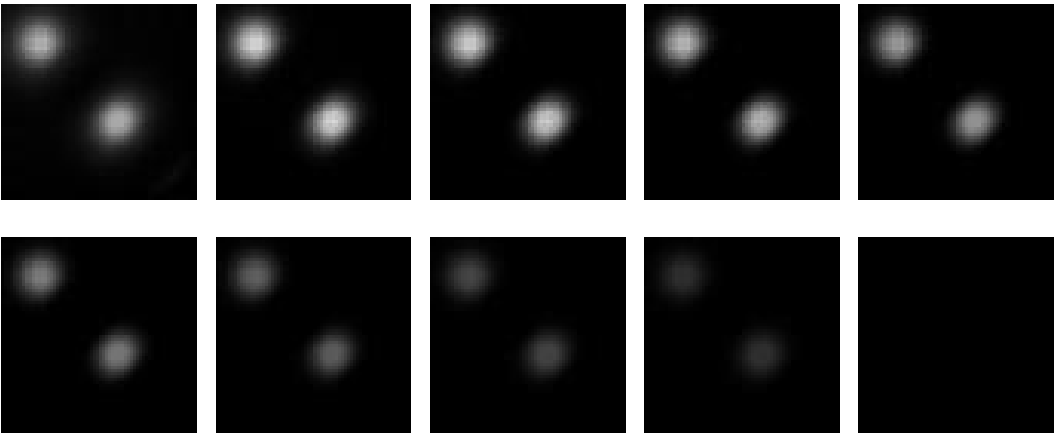

  \foreach \t in {1,2,3,4,5}{
  \includegraphics[width=0.17\linewidth]{figures//l_s_\t.jpg}
  }\\ \\
  \foreach \t in {6,7,8,9,10}{
  \includegraphics[width=0.17\linewidth]{figures//l_s_\t.jpg}
  }
  \caption{Source over time}
\end{figure}

\subsection{Unbalanced vector OMT}

We tested our approach on color image data. While the interpolations of density are color, the source layer is still gray scale.

\begin{figure}[H]
  \centering
  \includegraphics[width=0.2\linewidth]{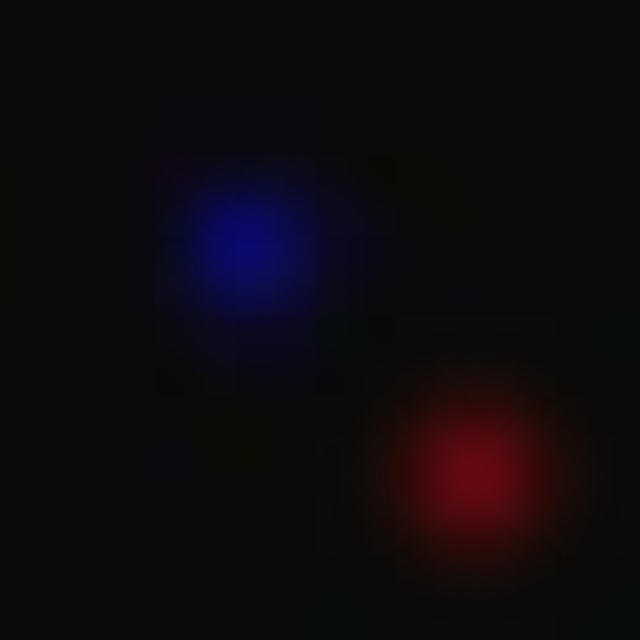}
  \includegraphics[width=0.2\linewidth]{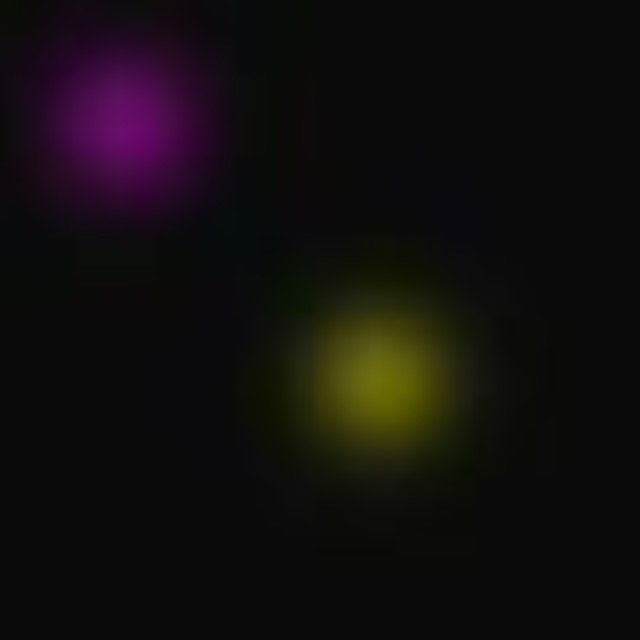}
  \caption{Source and target distributions}
\end{figure}
\begin{figure}[H]
  \foreach \t in {1,2,3,4,5}{
  \includegraphics[width=0.17\linewidth]{figures//VV4_\t.jpg}
  }\\ \\
  \foreach \t in {6,7,8,9,10}{
  \includegraphics[width=0.17\linewidth]{figures//VV4_\t.jpg}
  }
  \caption{Density over time}
\end{figure}
\begin{figure}[H]
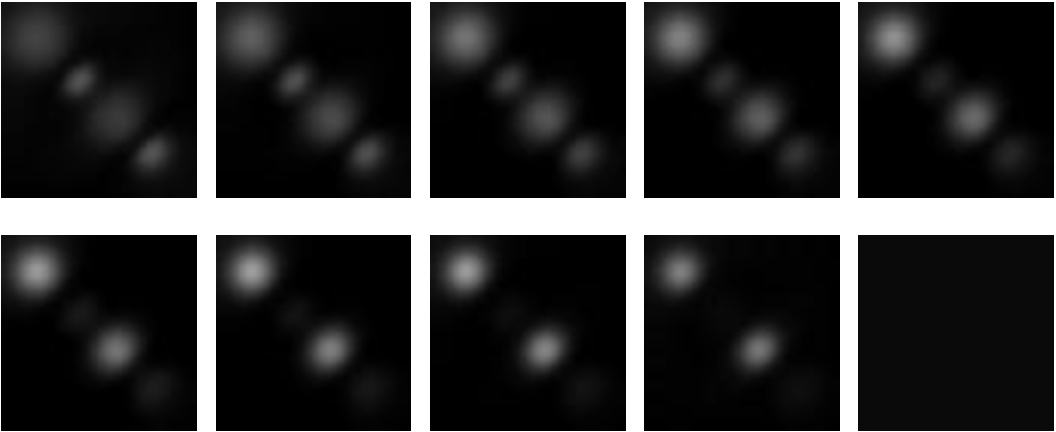

  \foreach \t in {1,2,3,4,5}{
  \includegraphics[width=0.17\linewidth]{figures//VV4_s_\t.jpg}
  }\\ \\
  \foreach \t in {6,7,8,9,10}{
  \includegraphics[width=0.17\linewidth]{figures//VV4_s_\t.jpg}
  }
  \caption{Source over time}
\end{figure}

\section{Conclusion}

Vector-valued OMT is a very powerful model. In this note, we reformulated unbalanced OMT and unbalanced vector-valued OMT by adding the source as a new layer. We gave a new way to include a source term and we proposed a very simple way to implement unbalanced OMT from general vector-valued OMT code. In the present work,  we only considered one kind of unbalanced formulation. Other flow-based unbalanced OMT settings fit into our model as well.  We believe our generalization of vector-valued OMT has not reached its full potential. Indeed, we only give different weight parameters to different layers. But what if we give different weight parameters to different areas, different nodes or different edges? We may use different parameters according to specific applications. This will be the subject of future research.

\section{Acknowledgments}

This study was supported by AFOSR grants (FA9550-17-1-0435, FA9550-20-1-0029), NIH grant (R01-AG048769), MSK Cancer Center Support Grant/Core Grant (P30 CA008748), and a grant from Breast Cancer Research Foundation (grant BCRF-17-193).

\bibliographystyle{plain}
\bibliography{refs}

\begin{thebibliography}{10}

\bibitem{Arjovsky2017}
Martin Arjovsky, Soumith Chintala, and L{\'{e}}on Bottou.
\newblock {Wasserstein GAN}.
\newblock 2017.

\bibitem{BB}
Jean-David Benamou and Yann Brenier.
\newblock A computational fluid mechanics solution to the {M}onge-{K}antorovich
  mass transfer problem.
\newblock {\em Numerische {M}athematik}, 84(3):375--393, 2000.

\bibitem{inbook}
Yongxin Chen, Tryphon Georgiou, and Allen Tannenbaum.
\newblock {\em Wasserstein Geometry of Quantum States and Optimal Transport of
  Matrix-Valued Measures}, pages 139--150.
\newblock 01 2018.

\bibitem{cheGeoTan16b}
Yongxin Chen, Tryphon~T Georgiou, and Allen Tannenbaum.
\newblock Interpolation of density matrices and matrix-valued measures: The
  unbalanced case.
\newblock {\em Euro. Jnl of Applied Mathematics}, 30(3):458--480, 2018.

\bibitem{CheGeoTan16}
Yongxin Chen, Tryphon~T Georgiou, and Allen Tannenbaum.
\newblock Matrix optimal mass transport: a quantum mechanical approach.
\newblock {\em IEEE Trans. Automatic Control}, 63(8):2612 -- 2619, 2018.

\bibitem{vectorvalued}
Yongxin Chen, Tryphon~T Georgiou, and Allen Tannenbaum.
\newblock Vector-valued optimal mass transport.
\newblock {\em SIAM Journal Applied Mathematics}, 78(3):1682--1696, 2018.

\bibitem{CheKaoHabGeoTan17}
Yongxin Chen, Kaoru Yamamoto, Eldad Haber, Tryphon~T. Georgiou, and Allen
  Tannenbaum.
\newblock An efficient algorithm for matrix-valued and vector-valued optimal
  mass transport.
\newblock {\em in preparation}, 2017.

\bibitem{chizat:hal-01271981}
Lenaic Chizat, Gabriel Peyr{\'e}, Bernhard Schmitzer, and Fran{\c c}ois-Xavier
  Vialard.
\newblock {Unbalanced Optimal Transport: Geometry and Kantorovich Formulation}.
\newblock working paper or preprint, August 2015.

\bibitem{Chizat}
Lenaic Chizat, Gabriel Peyr\'{e}, Bernhard Schmitzer, and Francois-Xavier
  Vialard.
\newblock An interpolating distance between optimal transport and
  {F}isher-{R}ao metrics.
\newblock {\em Foundations of Computational Mathematics}, 10:1--44, 2016.

\bibitem{Gangbo2019}
Wilfrid Gangbo, Wuchen Li, Stanley Osher, and Michael Puthawala.
\newblock {Unnormalized optimal transport}.
\newblock {\em Journal of Computational Physics}, 399:1--19, 2019.

\bibitem{Haker2004}
Steven Haker, Lei Zhu, Allen Tannenbaum, and Sigurd Angenent.
\newblock {Optimal mass transport for registration and warping}.
\newblock {\em International Journal of Computer Vision}, 60(3):225--240, 2004.

\bibitem{Statement2020}
James~C. Mathews, Saad Nadeem, Maryam Pouryahya, Zehor Belkhatir, Joseph~O.
  Deasy, Arnold~J. Levine, and Allen~R. Tannenbaum.
\newblock {Functional network analysis reveals an immune tolerance mechanism in
  cancer}.
\newblock {\em Proceedings of the National Academy of Sciences},
  117(28):16339--16345, jul 2020.

\bibitem{Otto}
Felix Otto.
\newblock The geometry of dissipative evolution equations: the porous medium
  equation.
\newblock {\em Communications in Partial Differential Equations}, 2001.

\bibitem{rr}
Svetlozar~T Rachev and Ludger R{\"u}schendorf.
\newblock {\em Mass Transportation Problems: {V}olumes {I} and {II}}.
\newblock Springer Science \& Business Media, 1998.

\bibitem{villani1}
C{\'e}dric Villani.
\newblock {\em Topics in Optimal Transportation}.
\newblock American Mathematical Soc., 2003.

\bibitem{villani2}
C{\'e}dric Villani.
\newblock {\em Optimal Transport: Old and New}, volume 338.
\newblock Springer Science \& Business Media, 2008.

\end{thebibliography}

\end{document}